%%%%%%%%%%%%%%%%%%%%%%%%%%%%%%%%%%%%%%%%%%%%%%%%%%%%%%%%%%%%%%%%%%%%%%%%%%%%
%%% Title:  Differentiable perturbation of unbounded operators
%%% Author: Andreas Kriegl, Peter W. Michor
%%% Remark: AmSTeX, 9 pages, final version December 8, 2002
%%%%%%%%%%%%%%%%%%%%%%%%%%%%%%%%%%%%%%%%%%%%%%%%%%%%%%%%%%%%%%%%%%%%%%%%%%%%
% TeX-NMB program applied to the file from 2002.9.27;14:40 on 2002.9.27; 14:40
%fmt=amstex
\nonstopmode
\input amstex
\input amsppt.sty   
\hsize 30pc
\vsize 47pc
\def\nmb#1#2{#2}         % used for renumbering, TeX should ignore.
\def\cit#1#2{\ifx#1!\cite{#2}\else#2\fi} %for citing references
             %= to table of content, invoked by kms-book.sty
\def\idx{}               % for producing index, invoked by kms-book.sty
\def\ign#1{}             %=ignore, invisible entry for the index only
\redefine\o{\circ}

\define\al{\alpha}
\define\be{\beta}
\define\ga{\gamma}
\define\de{\delta}
\define\ep{\varepsilon}

\define\la{\lambda}
\define\rh{\rho}
\define\si{\sigma}

\define\De{\Delta}

\define\Om{\Omega}
\redefine\i{^{-1}}
\define\x{\times}
\define\Lip{\operatorname{\Cal Lip}}
\def\today{\ifcase\month\or
 January\or February\or March\or April\or May\or June\or
 July\or August\or September\or October\or November\or December\fi
 \space\number\day, \number\year}
\topmatter
\title Differentiable perturbation of unbounded operators
\endtitle
\author 
Andreas Kriegl, Peter W\. Michor  
\endauthor
\leftheadtext{A\. Kriegl, P\.W\. Michor}
\address
A\. Kriegl: Institut f\"ur Mathematik, Universit\"at Wien,
Strudlhofgasse 4, A-1090 Wien, Austria
\endaddress
\email Andreas.Kriegl\@univie.ac.at\endemail
\address
P\. W\. Michor: Institut f\"ur Mathematik, Universit\"at Wien,
Strudlhofgasse 4, A-1090 Wien, Austria; {\it and}:
Erwin Schr\"odinger Institute of Mathematical Physics, Boltzmanngasse 
9, A-1090 Wien, Austria.
\endaddress
\email Peter.Michor\@esi.ac.at \endemail
%\date {\today} \enddate
\thanks 
PWM was supported by FWF, Projekt P~14195~MAT. 
\endthanks
\keywords Perturbation theory, differentiable choice of eigenvalues 
\endkeywords
%\subjclassyear{2000} % for amstex 2.2
\subjclass\nofrills{\rm 2000}
 {\it Mathematics Subject Classification}.\usualspace
 Primary 26C10.\endsubjclass
%\subjclass 47A55, 47A75 \endsubjclass
\abstract If $A(t)$ is a $C^{1,\al}$-curve of unbounded self-adjoint 
operators with compact resolvents and common domain 
of definition, then the eigenvalues can be parameterized $C^1$ in $t$. 
If $A$ is $C^\infty$ then the eigenvalues can be parameterized twice 
differentiably.
\endabstract
\endtopmatter

\document

\proclaim{Theorem}
Let $t\mapsto A(t)$ for $t\in \Bbb R$ be a curve of unbounded self-adjoint 
operators in a Hilbert space with common domain of definition and with 
compact resolvent.

\roster
\item"(A)" If $A(t)$ is real analytic in $t\in \Bbb R$, then 
     the eigenvalues and the eigenvectors of $A(t)$ may be 
     parameterized real analytically in $t$. 
\item"(B)" If $A(t)$ is $C^\infty$ in $t\in \Bbb R$ and if
     no two unequal continuously parameterized  
     eigenvalues meet of infinite order at any $t\in \Bbb R$, 
                 then the eigenvalues and the eigenvectors can be 
     parameterized smoothly in $t$, on the whole parameter domain.
\item"(C)" If $A$ is $C^\infty$, then 
     the eigenvalues of $A(t)$ may be parameterized
     twice differentiably in $t$. 
\item"(D)" If $A(t)$ is $C^{1,\al}$ for some $\al>0$ in 
       $t\in \Bbb R$, then the eigenvalues of $A(t)$ may be 
       parameterized in a $C^1$ way in $t$. 
\endroster
\endproclaim

Part \therosteritem{A} is due to Rellich \cit!{10} in 1940, 
see also \cit!{2} and \cit!{6}, VII,~3.9.
Part \therosteritem{B} has been proved in \cit!{1},~7.8, see 
also \cit!{8},~50.16, in 1997; there we gave also a different proof 
of \therosteritem{A}. The purpose of this paper is to prove 
parts \therosteritem{C} and 
\therosteritem{D}. 

Both results cannot be improved to obtain a $C^{1,\be}$-parameterization of the
eigenvalues for some $\be>0$, by the first example
below. In our proof of \therosteritem{D} 
the assumption $C^{1,\al}$ cannot be weakened to 
$C^1$, see the second example. For finite dimensional Hilbert spaces 
part \therosteritem{D} has been proved under the assumption of 
$C^1$  by Rellich \cit!{11}, with a small inaccuracy
in the auxiliary theorem on p\.~48: Condition (4) must be more 
restrictive, otherwise the induction argument on p\.~50 is not valid, 
since the proof on p\.~52 relies on the 
fact that all values coincide at the point in question. 
A proof can also be found in \cit!{6}, II,~6.8.
We need a strengthened version of this result, 
thus our proof covers it also.

We thank T\. and M\. Hoffmann-Ostenhof and T\. Kappeler for their 
interest and hints. 

\subhead Definitions and remarks \endsubhead
That $A(t)$ is a real analytic, $C^\infty$, or $C^{k,\al}$  curve of 
unbounded operators means the following:
There is a dense subspace $V$ of the Hilbert space $H$ 
such that $V$ is the domain of definition of each $A(t)$, and such 
that $A(t)^*=A(t)$. 
Moreover, we require 
that $t\mapsto \langle A(t)u,v\rangle$ 
is real analytic, $C^\infty$, or $C^{k,\al}$ for each $u\in V$ and $v\in H$. 
This implies that $t\mapsto A(t)u$ is of the same class $\Bbb R\to H$ 
for each $u\in V$ by \cit!{8},~2.3 or \cit!{5},~2.6.2. 
This is true because $C^{k,\al}$ can be described by 
boundedness conditions only; and for these the uniform boundedness 
principle is valid.  
A function $f$ is called $C^{k,\al}$ if it is $k$ times differentiable 
and for the $k$-th
derivative the expression $\frac{f^{(k)}(t)-f^{(k)}(s)}{|t-s|^\al}$ is 
locally bounded in $t\ne s$.

A sequence of continuous, real analytic, smooth, or twice 
differentiable functions 
$\la_i$ is said to {\it parameterize the eigenvalues, if
for each $z\in \Bbb R$ 
the cardinality $|\{i: \la_i(t)=z\}|$ equals
the multiplicity of $z$ as eigenvalue of $A(t)$.} 

The proof will moreover furnish the following (stronger) versions:
{\it\roster
\item"(C1)" If $A(t)$ is 
       $C^{3n,\al}$ in $t$ and if the multiplicity of an eigenvalue 
       never exceeds $n$, then the eigenvalues of 
       $A$ may be parameterized twice differentiably.
\item"(C2)"  If the multiplicity of 
       any eigenvalue never exceeds $n$, and if the resolvent 
       $(A(t)-z)\i$ is $C^{3n}$ into  
       $L(H,H)$ in $t$ and $z$ jointly, then the eigenvalues of 
       $A(t)$ may be parameterized twice differentiably in 
       $t$.  
\item"(D1)"  If the resolvent $(A(t)-z)\i$ is $C^{1}$ into 
       $L(H,H)$ in $t$ and $z$ jointly, then the eigenvalues of 
       $A(t)$ may be parameterized in a $C^1$ way in 
       $t$.  
\item"(D2)"  In the situations of 
       \therosteritem{D} and \therosteritem{D1} the 
       following holds: For any continuous parameterization $\la_i(t)$ 
       of all eigenvalues of $A(t)$, each function $\la_i$ has right 
       sided derivative $\la^{(+)}_i(t)$ and left sided one 
       $\la^{(-)}_i(t)$ at each $t$, and 
       $\{\la^{(+)}_i(t):\la_i(t)=z\}$ equals 
       $\{\la^{(-)}_i(t):\la_i(t)=z\}$ with correct multiplicities.
\endroster}

\subhead   Open problem \endsubhead
Construct a $C^1$-curve of unbounded self-adjoint operators with 
common domain and compact resolvent such that the eigenvalues cannot 
be arranged $C^1$. 

\subhead   Applications \endsubhead
Let $M$ be a compact manifold and let $t\mapsto g_t$ be a 
smooth curve of smooth Riemannian metrics on $M$. Then we get 
the corresponding smooth curve $t\mapsto \De(g_t)$
of Laplace-Beltrami operators on $L^2(M)$. By theorem (C) the 
eigenvalues can be arranged twice differentiably. 

Let $\Om$ be a bounded region in $\Bbb R^n$ with smooth boundary, and 
let $H(t)=-\De + V(t)$ be a $C^{1,\al}$-curve of Schr\"odinger 
operators with varying potential and Dirichlet boundary conditions.
Then the eigenvalues can be arranged $C^1$. 

\subhead Example
\endsubhead 
This is an elaboration of \cit!{1},~7.4.
Let $S(2)$ be the vector space of all symmetric real 
$(2\x 2)$-matrices. 
We use the general curve lemma \cit!{8},~12.2:
{\it
There exists a converging sequence of reals $t_n$ with the following 
property:
Let $A_n\in C^\infty(\Bbb R,S(2))$ be any sequence of functions which 
converges  
fast to 0, i.e., for each $k\in \Bbb N$ the sequence $n^kA_n$ is 
bounded in $C^\infty(\Bbb R,S(2))$. 
Then there exists a smooth curve $A\in C^\infty(\Bbb R,S(2))$ 
such that $A(t_n+s)=A_n(s)$ for 
$|s|\le \frac1{n^2}$, for all $n$.}

We use it for 
$$
A_n(t) := \pmatrix
        \frac 1{2^{n^2}} & \frac t{2^n} \\
        \frac t{2^n} & -\frac 1{2^{n^2}}\\
        \endpmatrix
= \frac 1{2^{n^2}} \pmatrix
        1 & \frac t{s_n} \\
        \frac t{s_n} & -1\\
        \endpmatrix, \text{ where } s_n := 2^{n-n^2}\le \frac 1{n^2}.
$$
The eigenvalues of $A_n(t)$ and their derivatives are
$$
\la_n(t) = \pm\frac 1{2^{n^2}} \sqrt{1+(\tfrac t{s_n})^2},\quad
\la_n'(t) = \pm\frac {2^{n^2-2n}t}{\sqrt{1+(\frac t{s_n})^2}}.
$$
Then
$$\align
\frac{\la'(t_n+s_n)-\la'(t_n)}{s_n^\al} &=
\frac{\la_n'(s_n)-\la_n'(0)}{s_n^\al}
=\pm\frac {2^{n^2-2n}s_n}{s_n^\al\sqrt{2}}\\
&=\pm\frac{2^{n(\al(n-1)-1)}}{\sqrt{2}}  \to \infty \text{ for }\al>0.
\endalign$$
By \cit!{1},~2.1, we may always find a 
twice differentiable square root of a non-negative smooth function, 
so that the eigenvalues $\la$ are functions 
which are twice differentiable but not $C^{1,\al}$ for any $\al>0$. 

Note that the normed eigenvectors cannot be chosen continuously in 
this example (see also example \cit!{9}, \S 2). Namely, we have
$$
A(t_n)=A_n(0)=\frac1{2^{n^2}}\pmatrix 1 & 0 \\ 0 &-1 \endpmatrix,\qquad
A(t_n+s_n)=A_n(s_n)=\frac1{2^{n^2}}\pmatrix 1 & 1 \\ 1 &-1 \endpmatrix.
$$

\proclaim{Resolvent Lemma}
If $A$ is $C^{k,\al}$ for some $1\le k\le \infty$ and $\al>0$, then the resolvent 
$(t,z)\mapsto (A(t)-z)\i\in L(H,H)$ is $C^{k}$ on its natural domain. 
\endproclaim

By $C^{\infty,\al}$ we mean $C^\infty$.

\demo{\bf Proof}
By definition the function $t\mapsto \langle A(t)v,u \rangle$ is of 
class $C^{k,\al}$ for each $v\in V$ and $u\in H$. Then by \cit!{4},~5 or 
\cit!{8},~2.3 (extended from $C^{k,1}$ to $C^{k,\al}$ with 
essentially the same proof),  
the curve $t\mapsto A(t)v$ is of class $C^{k,\al}$ into $H$.

For each $t$ consider the norm $\|u\|_t^2:=\|u\|^2+\|A(t)u\|^2$ on 
$V$.
Since $A(t)=A(t)^*$ is closed, $(V,\|\quad\|_t)$ is again a 
Hilbert space with inner product 
$\langle u,v\rangle_t:=\langle u,v\rangle+\langle A(t)u,A(t)v\rangle$. 

{\it \therosteritem{\nmb.{1}} Claim. 
All these norms $\|\quad\|_t$ on $V$ are equivalent, 
locally uniformly in $t$.} 
We then equip $V$ with one of the 
equivalent Hilbert norms, say $\|\quad\|_0$.

Note first that $A(t):(V,\|\quad\|_s)\to H$ is bounded since the 
graph of $A(t)$ is closed in $H\x H$, contained in $V\x H$ and thus 
also closed in $(V,\|\quad\|_s)\x H$. 
For fixed $u,v\in V$, the function
$t\mapsto \langle u,v\rangle_t=\langle u,v \rangle+\langle A(t)u,A(t)v 
\rangle$ is $C^{k,\al}$ since $t\mapsto A(t)u$ is it.
Thus it is also locally Lipschitz ($C^{0,1}=\Lip^0$).
By the multilinear uniform boundedness principle 
(\cit!{8},~5.18) or \cit!{5},~3.7.4) the 
mapping $t\mapsto \langle \quad,\quad\rangle_t$ is $C^{0,1}$ into the 
space of bounded bilinear forms on $(V,\|\quad\|_s)$ for each fixed 
$s$. By the exponential law \cit!{5},~4.3.5 for $\Lip^0$ the mapping 
$(t,u,v)\mapsto \langle u,v\rangle_t$ is $C^{0,1}$ 
from $\Bbb R\x (V,\|\quad\|_s)\x (V,\|\quad\|_s)\to \Bbb R$ for each 
fixed $s$. 
Therefore and by homogeneity in $(u,v)$ the set 
$\{\|u\|_t:|t|\le K,\|u\|_s\le 1 \}$ is bounded by 
some $L_{K,s}$ in $\Bbb R$.
Thus $\|u\|_t\le L_{K,s}\|u\|_s$ for all $|t|\le K$, i.e\. 
all Hilbert norms $\|\quad\|_t$ are locally uniformly 
equivalent, and claim \therosteritem{\nmb!{1}} follows. 

By \cit!{4},~5 and the  
linear uniform boundedness theorem 
we see that $t\mapsto A(t)$ is a $C^{k,\al}$-mapping $\Bbb R\to L(V,H)$, and 
thus is $C^{k}$ in the usual sense, again by \cit!{4},~5. 
Alternatively, if reference \cit!{4} is not available, one may use
\cit!{8},~2.3, extended from $C^{k,1}$ to $C^{k,\al}$ with 
essentially the same proof, and note that it suffices to test with linear 
mappings which recognize bounded sets, 
by \cit!{8},~5.18. 
Alternatively again, one may use  
\cit!{5},~3.7.4~+~4.1.12, extended from $C^{k,1}$ to $C^{k,\al}$.

If for some $(t,z)\in \Bbb R\x\Bbb C$ the bounded operator 
$A(t)-z:V\to H$ is invertible, then this is true locally and 
$(t,z)\mapsto (A(t)-z)\i:H\to V$ is $C^{k}$, by the chain rule,
since inversion is smooth on the Banach space $L(V,H)$.
\qed\enddemo

Since each $A(t)$ is Hermitian with compact resolvent 
the \idx{\it global resolvent set}
$\{(t,z)\in \Bbb R\x\Bbb C: (A(t)-z):V\to H \text{ is invertible}\}$  
is open and connected.
Moreover, $(A(t)-z)\i:H\to H$ is a compact operator 
for some (equivalently any) $(t,z)$ if and only if the inclusion 
$i:V\to H$ is compact, since $i=(A(t)-z)\i\o(A(t)-z): V\to H\to H$.

\subhead Resolvent example \endsubhead 
The resolvent lemma cannot be 
improved. We describe a curve $A(t)$ of self adjoint unbounded operators 
on $\ell^2$ with 
compact resolvent and common domain $V$ of definition, such that 
$t\mapsto \langle A(t)v,u \rangle$ is $C^1$ for all $v\in V$ and 
$u\in \ell^2$, but $t\mapsto A(t)$ is even not differentiable at 0 into 
$L(V,\ell^2)$.  

Let $\la_1\in C^\infty(\Bbb R,\Bbb R)$ be nonnegative with 
compact support and $\la_1'(0)=0$. We consider the multiplication 
operator $B(t)$ on $\ell^2$ given on the standard basis $e_n$ by 
$B(t)e_n:=(1+\frac1n\la_1(nt))e_n=:\la_n(t)e_n$ which is bounded with 
bounded inverse. Then the function  $t\mapsto \langle B(t)x,y\rangle$ 
is $C^1$ with derivative $\langle B_1(t)x,y \rangle$, where $B_1(t)$ 
is given by $B_1(t)e_n=\la_n'(t)e_n=\la_1'(nt)e_n$, since for fixed 
$t$ we have that
$$
\mu_n(s):=\frac{\la_n(t+s)-\la_n(t)}{s}-\la_n'(t)
=\frac{\la_1(nt+ns)-\la_1(nt)}{ns}-\la_1'(nt)
$$
converges to 0 for 
$s\to 0$ pointwise in $n$ and is bounded uniformly in $n$:
$$
\Bigl|\sum_n\mu_n(s)x_ny_n\Bigr| 
\le \sum_{n=N+1}^\infty|\mu_n(s)x_ny_n| 
     + \sum_{n=1}^N|\mu_n(s)x_ny_n| 
\le \sup_n|\mu_n(s)|\ep + \ep \|x\|\|y\|.
$$
Moreover, $\langle B_1(t)x,y \rangle$ is continuous in $t$ since 
$\la_n'(t+s)-\la_n'(t)=\la_1'(nt+ns)-\la_1'(nt)$ 
also converges to 0 for 
$s\to 0$ pointwise in $n$ and is bounded uniformly in $n$.
But $t\mapsto B(t)$ is not differentiable at $0$ into 
$L(\ell^2,\ell^2)$ since 
$$
\Bigl\|\frac{B(t)-B(0)}{t}-B_1(0)\Bigr\| 
=\sup_n\Bigl|\frac{\la_n(t)-\la_n(0)}{t}-\la_n'(0)\Bigr|
=\sup_n\Bigl|\frac{\la_1(nt)-\la_1(0)}{nt}-\la_1'(0)\Bigr|
$$
is bounded away from 0, for $t\to 0$.
Finally, let $C:\ell^2\to \ell^2$ be the compact invertible given by
$Ce_n=\frac1n e_n$.
We take $V=C(\ell^2)$, and $A(t)=B(t)\o C\i$. 

\demo{\bf Proof of the theorem}
By the resolvent lemma, 
\therosteritem{D1} implies \therosteritem{D}, and likewise  
\therosteritem{C2} implies 
\therosteritem{C} and \therosteritem{C1}. 
\enddemo

\demo{\bf Proof of \therosteritem{D1}} 

(\nmb.{2}) {\it Claim. If $f:(a,b)\to \Bbb R$ is continuous 
and if $f(t)$ has only finitely many cluster points for $t\to b$ 
then the limit $\lim_{t\nearrow b}f(t)$ exists.} 
Otherwise, by the intermediate value theorem, we have a whole 
interval of cluster points. 

{\it \therosteritem{\nmb.{3}} Claim. Let $z$ be an eigenvalue of 
$A(s)$ of multiplicity $N$. Then there exists an open box 
$(s-\de,s+\de)\x(z-\ep,z+\ep)$ and $C^1$-functions 
$\mu_1,\dots,\mu_N: (s-\de,s+\de)\to (z-\ep,z+\ep)$ which 
parameterize all eigenvalues $\la$ with $|\la-z|<\ep$ of $A(t)$ 
for $|t-s|<\de$ with correct multiplicities.} 

We choose a simple closed smooth curve $\ga$ in the resolvent set of 
$A(s)$ for fixed $s$ enclosing only $z$ among all eigenvalues of $A(s)$. 
Since the global resolvent set is open, no eigenvalue 
of $A(t)$ lies on $\ga$, for $t$ near $s$.
Since 
$$
t\mapsto -\frac1{2\pi i}\int_\ga (A(t)-z)\i\;dz =: P(t,\ga) = P(t)
$$
is a $C^1$ curve of projections (on the direct sum of all 
eigenspaces corresponding to eigenvalues in the interior of $\ga$) 
with finite dimensional ranges, the ranks
(i\.e\. dimension of the ranges) must be constant: it is easy to see 
that the  
(finite) rank cannot fall locally, and it cannot increase, since the 
distance in $L(H,H)$ of $P(t)$ to the subset of operators of 
rank $\le N=\operatorname{rank}(P(s))$ is continuous in $t$ and is either 
0 or 1.   
So for $t$ near $s$, say $t\in I:=(s-\de,s+\de)$, 
there are equally many eigenvalues in the 
interior of $\ga$, and we may call them 
$\la_i(t)$ for $1\le i\le N$ (repeated 
with multiplicity), so that each $\la_i$ is continuous 
(this is well known and follows easily from the proof of
\therosteritem{C2}).

Then the image of $t\mapsto P(t,\ga)$, for $t$ 
near $s$, describes a $C^1$ finite dimensional vector subbundle of 
$\Bbb R\x H\to \Bbb R$, since its rank is constant.
For each $t$ choose an 
orthonormal system of eigenvectors $v_j(t)$ of $A(t)$ corresponding to 
these $\la_j(t)$.
They form a (not necessarily continuous) framing of 
this bundle. 
For any $t$ near $s$ and any sequence $t_k\to t$ there is a 
subsequence again denoted by $t_k$ such that each 
$v_j(t_k)\to w_j(t)$ where the $w_i(t)$ form again an orthonormal system of 
eigenvectors of $A(t)$ for the sum $P(t)(H)$ of the eigenspaces of 
the $\la_i(t)$ (Here we use the local triviality of the vector bundle).
Now consider 
$$
\frac{A(t)-\la_i(t)}{t_k-t}v_i(t_k) +
\frac{A(t_k)-A(t)}{t_k-t}v_i(t_k) -
\frac{\la_i(t_k)-\la_i(t)}{t_k-t}v_i(t_k) = 0.
\tag{\nmb.{4}}$$
For $t=s$ we
take the inner product of \therosteritem{\nmb!{4}} with each
$w_j(s)$, note that then the  
first summand vanishes since all $\la_i(s)$ agree, and let $k\to \infty$ 
to obtain that (for $j\ne i$) the $w_i(s)$ are a basis of eigenvectors of 
$P(s)A'(s)|P(s)(H)$ with eigenvalues (for $j=i$) 
$\lim_{k}\frac{\la_i(t_k)-\la_i(s)}{t_k-s}$. 
By \therosteritem{\nmb!{2}}, 
$$
\lim_{h\searrow 0}\frac{\la_i(s+h)-\la_i(s)}{h} = \rh_i,
$$
where the $\rh_i$ are the eigenvalues of $P(s)A'(s)|P(s)(H)$ (with correct
multiplicities).
So the right handed derivative $\la^{(+)}_j(s)$ of each $\la_j$ exists 
at $s$. Similarly the left handed derivative $\la^{(-)}_j(s)$ 
exists, 
and they form the same set of numbers with the correct multiplicities. 
Thus there exists a permutation $\si$ of $\{1,\dots,N\}$ such that 
the
$$
\nu_i(t) := \cases \la_i(t) &\text{ for }t\le s \\
                  \la_{\si(i)}(t) &\text{ for }t\ge s \endcases
\tag{\nmb.{5}}$$
{\it parameterize all eigenvalues in the box by continuous 
functions which are differentiable at $s$.}  

For $t\ne s$, take the  
inner product of \therosteritem{\nmb!{4}} with $w_i(t)$ to conclude that
\roster
\item"\thetag{\nmb.{6}}"
{\it 
$\la_i^{(+)}(t) = \langle A'(t)w_i(t),w_i(t) \rangle$
for a unit eigenvector $w_i(t)$ of $A(t)$
with eigenvalue $\la_i(t)$}.
\endroster

Now we show claim \therosteritem{\nmb!{3}} by induction on $N$. 
Let $t_1\in I$ be such that not all $\la_i(t_1)$ agree. 
Then $\{1,\dots,N\}$ decomposes into the subsets $\{i:\la_i(t_1)=w\}$.
Then for $i$ and $k$ in different subsets 
$\la_i(t)\ne\la_k(t)$ for all $t$ in an open interval $I_1$ 
containing $t_1$. Thus by induction on each subset  
\therosteritem{\nmb!{3}} holds on $I_1$. 

Next let $I_2\subseteq I$ be an open interval containing only points 
$t_1$ as above. Let $J$ be a maximal open subinterval on which 
\therosteritem{\nmb!{3}} holds. Assume for contradiction that 
the right (say) endpoint $b$ of 
$J$ belongs to $I_2$, then there is a $C^1$-parameterization of all 
$N$ eigenvalues on an open interval $I_b$ containing $b$ by the argument 
above. Let $t_2\in J\cap I_b$. Renumbering the $C^1$ parameterization 
to the right of $t_2$ suitably we may extend the $C^1$ 
parameterization beyond $b$, a contradiction. Thus 
\therosteritem{\nmb!{3}} holds on $I_2$. 

Now we consider the closed set $E=\{t\in I: 
\la_1(t)=\dots=\la_N(t)\}$. Then $I\setminus E$ is open, thus a 
disjoint union of open intervals on which there exists a 
$C^1$-parameterization $\mu_i$ of all eigenvalues. 
Consider first the set $E'$ of all isolated points in $E$. Then 
$E'\cup (I\setminus E)$ is again open and thus a disjoint union of open 
intervals, and for each point $t\in E'$ we apply in turn the following
arguments: extending all $\mu_i$'s by the single value at $t$ we get a
continuous extension near $t$. Then by \therosteritem{\nmb!{5}}, 
we may renumber the $\mu_i$ 
to the right of $t$ in such a way that they fit together differentiably at 
$t$. The derivatives are also continuous at $t$: They have only finitely
many clusterpoints for $t_k\to t$ by applying \therosteritem{\nmb|{6}} to
$t_k$ and choosing a subsequence such that the $w_i(t_k)$ converge. 
Now we apply the arguments surrounding \therosteritem{\nmb|{4}} with the
$v_j(t_k)$ replaced by $w_j(t_k)$ to conclude that
\therosteritem{\nmb|{6}} converges to $\rh_i(t)=\mu_i'(t)$.
Thus \therosteritem{\nmb!{3}} holds on 
$E'\cup (I\setminus E)$. 

We extend each 
$\mu_i$ to the whole of $I$ by taking the single continuous function 
on $E\setminus E'$. Let $t\in E\setminus E'$. Then for the 
parameterization $\nu_i$ of \therosteritem{\nmb!{5}} of all 
eigenvalues which is differentiable at $t$ all derivatives 
$\nu_i'(t)$ agree since $t$ is a cluster point of $E$. 
Thus also $\mu_i'(t)$ exists and equals $\nu_i'(t)$. 
So all $\mu_i$ are differentiable on $I$. 

To see that $\mu_i'$ is continuous at $t\in E\setminus E'$, let 
$t_n\to t$ be such that $\mu_i'(t_n)$ converges (to a cluster point 
or $\pm\infty$). 
Then by \therosteritem{\nmb!{6}} we have 
$\mu_i'(t_k)=\langle A'(t_k)w_i(t_k),w_i(t_k) \rangle$ for 
eigenvectors $w_i(t_k)$ of $A(t_k)$ with eigenvalue $\mu_i(t_k)$. 
Passing to a subsequence we may assume that the $w_i(t_k)$ converge 
to an orthonormal basis of eigenvectors of $A(t)$, then 
$\langle A'(t_k)w_i(t_k),w_i(t_k) \rangle$ converges to some of the 
equal eigenvalues $\rh_i$ of $P(t)A'(t)|P(t)(H)$  which also equal 
the $\nu_i'(t)$. 

So \therosteritem{\nmb!{3}} is completely proved.

{\it \therosteritem{\nmb.{7}} Claim. 
Let $I$ be a compact interval.
Let $t\mapsto\la_i(t)$ be a differentiable eigenvalue of $A(t)$, 
defined on some subinterval of $I$.
Then
$$
|\la_i(t_1)-\la_i(t_2)| \le (1+|\la_i(t_2)|)(e^{a|t_1-t_2|}-1) 
$$
holds for a positive constant $a$ depending only on $I$.}

From \thetag{\nmb!{6}} we conclude, where $V_t=(V,\|\quad\|_t)$, 
$$\align
|\la_i'(t)|&\le \|A'(t)\|_{L(V_t,H)}
     \|w_i(t)\|_{V_t}\|w_i(t)\|_H \\
&= \|A'(t)\|_{L(V_t,H)}
     \sqrt{\|w_i(t)\|_H^2+\|A(t)w_i(t)\|_H^2}\cdot 1 \\
&= \|A'(t)\|_{L(V_t,H)}\sqrt{1+\la_i(t)^2}\le C+C|\la_i(t)|,\\
\endalign$$
for a constant $C$ 
since all norms $\|\quad\|_t$ are locally in $t$ uniformly 
equivalent, see claim \therosteritem{\nmb!{1}} above.
By Gronwall's lemma (see e\.g\. \cit!{3}, (10.5.1.3)) this implies 
claim \therosteritem{\nmb!{7}}.
%Namely,
%$$\align
%|\la(t)|-|\la(0)|&\le |\la(t)-\la(0)|\le \int_0^t |\la'(s)|ds 
%\le \int_0^t C(1+|\la(s)|)ds, \\
%|\la(t)| &\le |\la(0)|+Ct+ \int_0^t C|\la(s)|ds, \\
%w(t) \le \ph(t) &+  \int_0^t \ps(s)w(s)ds \Rightarrow 
%w(t) \le \ph(t)+\int_0^t\ph(s)\ps(s)e^{\int_s^t \ps }ds\\
%|\la(t)| &\le |\la(0)|+Ct+ \int_0^t (|\la(0)|+Cs)Ce^{C(t-s)}ds\\
%&= |\la(0)|+Ct -1 -Ct -|\la(0)| + e^{Ct}(1+|\la(0)|) \\
%&= |\la(0)| + (e^{Ct}-1)(1+|\la(0)|). 
%\endalign$$

By the following arguments we can conclude that all eigenvalues may 
be parameterized in a $C^1$ way.
Let us first number all eigenvalues of $A(0)$ (increasingly, say).

We consider families of $C^1$-functions $(\mu_i)_{i\in \al}$ indexed by 
ordinals $\al$,
defined on open intervals $I_i$ containing some fixed $t_0$,
which parameterize eigenvalues.

The set of all these sequences is partially ordered by inclusion of ordinals
and then by restriction of the component functions. 
Obviously for each increasing chain of such 
sequences the union is again such a sequence. By Zorn's lemma there 
exists a maximal family $(\mu_i)$.

We claim that for any maximal family each component function $\mu_i$ 
is globally defined: If  
not let $b<\infty$ be the (right, say) boundary point of $I_i$. By 
claim \therosteritem{\nmb!{7}} the limit $\lim_{t\nearrow b}\mu_i(t)=:z$ 
exists. By claim \therosteritem{\nmb!{3}} there exists a box 
$(b-\de,b+\de)\x (z-\ep,z+\ep)$ such that all eigenvalues $\la$ with 
$|\la-z|<\ep$ of $A(t)$ for $|t-b|<\de$ are parameterized by 
$C^1$ functions $\la_i:(b-\de,b+\de)\to (z-\ep,z+\ep)$ (with 
multiplicity). 
Consider the $\mu_j$ hitting this box (at the vertical boundaries only).
The endpoints of the corresponding
intervals $I_j$ give a partition of $(b-\de,b+\de)$ into finitely many 
subintervals. 
We apply the lemma below on each subinterval and glue at the ends of the
subintervals in
$C^1$-fashion using \therosteritem{\nmb|{5}} to 
obtain an extension of at least $\mu_i$, so the 
family was not maximal. 

Finally we claim that any maximal family $(\mu_i)$ parameterizes 
all eigenvalues of $A(t)$ with right multiplicities, for each $t\in \Bbb R$. 
If not, there is an eigenvalue $z$ of $A(t_0)$ with 
$|\{i: \mu_i(t_0)=z\}|$ less than the multiplicity of $z$. By claim 
\therosteritem{\nmb!{3}} and the lemma below 
we can then conclude 
again that the sequence was not maximal. 
\qed\enddemo

\proclaim{Lemma}
Suppose that $\la_{1},\dots,\la_N$ are
real-valued $C^1$ (twice differentiable) 
functions defined on an interval $I$, and 
that $\mu_1,\dots,\mu_k$ for $k\le N$ are also $C^1$ (twice differentiable) 
functions on $I$ such 
that $|\{j:\mu_j(t)=z\}|\le |\{i:\la_i(t)=z\}|$ for all $t\in I$ and 
$z\in \Bbb R$. 

Then there exist $C^1$ (twice differentiable) 
functions $\mu_{k+1},\dots,\mu_N$ on $I$ 
such that for all $t\in I$ and $z\in \Bbb R$ we have 
$|\{j:1\le j\le N,\mu_j(t)=z\}| = |\{i:\la_i(t)=z\}|$ . 
\endproclaim

\demo{\bf Proof} 
We treat the case $C^1$ and indicate the necessary changes 
in brackets for the twice
differentiable case.

We use induction on $N$. 
Let us assume that 
the statement is true if the number of functions is less than $N$. 

First suppose that for given $t_1\in I$ not all $\la_i(t_1)$ agree. 
Then for $i\in \{k:\la_1(t_1)=\la_k(t_1)\}\not\ni j$ we 
have $\la_i(t)\ne \la_j(t)$ for all $t$ in an open interval $I_1$ 
containing $t_1$, and similarly for the $\mu_j$. 
Thus by induction for both groups the statement 
holds on $I_1$.   

Now suppose that for no point $t$ in $I$ we have 
$\la_1(t)=\dots=\la_N(t)$. Let $I_1$ be a maximal open subinterval of $I$ 
for which the statement is true with functions $\mu^1_i$ for $i>k$. 
Assume for contradiction that the 
right (say) endpoint $b$ of $I_1$ is an interior point of $I$. 
By the first case,
the statement holds for an open neighborhood $I_2$ of $b$, with 
functions $\mu^2_i$ for $i>k$. 
Let $t_0\in I_1\cap I_2$. 
We may continue each solution $\mu_i^1$ in 
$\{t\in I_1:t\le t_0\}$ by a suitable solution $\mu_{\pi(i)}^2$ on 
$\{t\in I_2:t\ge t_0\}$ for a 
suitable permutation $\pi$: 
Let $t_m\nearrow t_0$. For every $m$ there exists a permutation
$\pi$ of $\{1,\dots,n\}$ such that 
$\mu^2_{\pi(i)}(t_m)=\mu^1_i(t_m)$ for all $i$.
By passing to a subsequence, again denoted $t_m$, we may assume
that the permutation does not depend on $m$.
By passing again to a subsequence we may also assume
that $(\mu^2_{\pi(i)})'(t_m)=(\mu^1_i)'(t_m)$ (and in the twice
differentiable case, again for a subsequence,
that $(\mu^2_{\pi(i)})''(t_m)=(\mu^1_i)''(t_m)$) for all $i$
and all $m$.
So we may paste $\mu^2_{\pi(i)}(t)$ for $t\ge t_0$
with $\mu^1_i(t)$ for $t<t_0$ to obtain a $C^1$ (twice 
differentiable) parameterization
on an interval larger than $I_1$, a contradiction.

In the general case,
we consider the closed set
$E=\{t\in I: \la_1(t)=\dots=\la_N(t)\}$.
Then $I\setminus E$ is open, thus a disjoint union of open intervals. 
By the second case the result holds on each of these 
open intervals.
Consider first the set $E'$ of all isolated points in $E$. Then 
$E'\cup (I\setminus E)$ is again open and thus a union of open 
intervals, and for each point $t\in E'$ we may renumber the $\mu_i$ 
to the right of $t$ in such a way that they fit together $C^1$ (twice
differentiable) at
$t$. Thus the result holds on 
$E'\cup (I\setminus E)$. 

We extend each 
$\mu_i$ to the whole of $I$ by taking the single continuous function 
on $E\setminus E'$. Let $t\in E\setminus E'$. Then 
all $\la_i'(t)=:\la'(t)$ agree since $t$ is cluster point of $E$ 
(and 
all $\la_i''(t)=:\la''(t)$ agree by considering second order 
difference quotients on points in $E$). 
Thus $\mu_i$ is (twice) differentiable at $t$ with $\mu_i'(t)=\la'(t)$ 
(and $\mu_i''(t)=\la''(t)$).

In the $C^1$ case we have still to check
that $\mu_i'$ is continuous at $t\in E\setminus E'$: Let 
$t_n\to t$, then 
$\mu_i'(t_n)=\la_{\si_n(i)}'(t_n)\to \la'(t)=\mu_i'(t)$. 
\qed\enddemo

\demo{\bf Proof of \therosteritem{C2}}
By assumption 
the resolvent $(A(t)-z)\i$ is
$C^{3n}$ jointly in $(t,z)$ where $n$ may be $\infty$.

{\it \therosteritem{\nmb.{8}} Claim. Let $z$ be an eigenvalue of 
$A(s)$ of multiplicity $N\le n$. Then there exists an open box 
$(z-\ep,z+\ep)\x(s-\de,s+\de)$ and twice differentiable functions 
$\mu_1,\dots,\mu_N: (s-\de,s+\de)\to (z-\ep,z+\ep)$ which 
parameterize all eigenvalues $\la$ with $|\la-z|<\ep$ of $A(t)$ 
for $|t-s|<\de$ with correct multiplicities.} 

We choose a simple closed smooth curve $\ga$ in the resolvent set of 
$A(s)$ for fixed $s$ enclosing only $z$ among all eigenvalues of $A(s)$. 
As in the proof of claim \therosteritem{\nmb!{3}} 
we see that 
$$
t\mapsto -\frac1{2\pi i}\int_\ga (A(t)-z)\i\;dz =: P(t,\ga) = P(t)
$$
is a $C^{3n}$ curve of projections 
with finite dimensional ranges of constant rank.

So for $t$ near $s$, there are equally many eigenvalues in the 
interior of $\ga$, and we may call them 
$\mu_i(t)$ for $1\le i\le N$ (repeated 
with multiplicity).
Let us denote by 
$e_i(t)$ for $1\le i\le N$ a corresponding system of eigenvectors of $A(t)$.
Then by the residue theorem we have 
$$
\sum_{i=1}^{N}\mu_i(t)^p e_i(t)\langle e_i(t),\quad \rangle
     = -\frac1{2\pi i}\int_\ga z^p(A(t)-z)\i\;dz 
$$
which is  $C^{3n}$ in $t$ near $s$, as a curve of operators in $L(H,H)$ 
of rank $N$.

{\it \therosteritem{\nmb.{9}} Claim. 
Let $t\mapsto T(t)\in L(H,H)$ be a  $C^{3n}$ curve of 
operators of rank $N$ in Hilbert space such that $T(0)T(0)(H)=T(0)(H)$.
Then $t\mapsto \operatorname{Trace}(T(t))$ is  $C^{3n}$ near $0$.}

This is claim 2 from \cit!{1},~7.8 for $C^{3n}$ instead of $C^\infty$.
\comment
Let $F:=T(0)(H)$.
Then $T(t)=(T_1(t),T_2(t)):H\to F\oplus F^\bot$ and 
the image of $T(t)$ is the space 
$$\align
T(t)(H) &= \{(T_1(t)(x),T_2(t)(x)):x\in H\}\\
&= \{(T_1(t)(x),T_2(t)(x)):x\in F\}\text{ for }t\text{ near }0\\
&= \{(y,S(t)(y)):y\in F\}, \text{ where }S(t):=T_2(t)\o (T_1(t)|F)\i.\\
\endalign$$
Note that $S(t):F\to F^\bot$ is  $C^{3n}$ in $t$ by finite 
dimensional inversion for $T_1(t)|F:F\to F$.
Now
$$\align
\operatorname{Trace}(T(t)) &= \operatorname{Trace}\left(  
     \pmatrix 1 & 0 \\ -S(t) & 1 \endpmatrix 
     \pmatrix T_1(t)|F & T_1(t)|F^\bot \\ 
          T_2(t)|F & T_2(t)|F^\bot\endpmatrix
     \pmatrix 1 & 0 \\ S(t) & 1 \endpmatrix
     \right)\\
&= \operatorname{Trace}\left(  
     \pmatrix T_1(t)|F & T_1(t)|F^\bot \\ 0 
          & -S(t)T_1(t)|F^\bot+T_2(t)|F^\bot \endpmatrix
     \pmatrix 1 & 0 \\ S(t) & 1 \endpmatrix\right)\\
&= \operatorname{Trace}\left(  
     \pmatrix T_1(t)|F & T_1(t)|F^\bot \\ 0 & 0 \endpmatrix
     \pmatrix 1 & 0 \\ S(t) & 1 \endpmatrix\right),
     \text{ since rank}=N\\
&= \operatorname{Trace}
     \pmatrix T_1(t)|F+(T_1(t)|F^\bot)S(t)  & T_1(t)|F^\bot \\ 
          0 & 0 \endpmatrix\\
&= \operatorname{Trace}\Bigl(T_1(t)|F+(T_1(t)|F^\bot)S(t):F\to F\Bigr),
\endalign$$ 
which is visibly  $C^{3n}$ since $F$ is finite dimensional. 
\endcomment
We conclude that the Newton polynomials
$$
s_p(t):=\sum_{i=1}^{N}\mu_i(t)^p =
-\frac1{2\pi i}\operatorname{Trace}\int_\ga z^p(A(t)-z)\i\;dz 
$$
are  $C^{3n}$ for $t$ near $s$.
Hence also the elementary symmetric polynomials 
$$
\si_p(t)=\sum_{i_1<\dots <i_N} \mu_{i_1}(t)\dots\mu_{i_p}(t)
$$
are $C^{3n}$, and thus  
$\{\mu_i(t): 1\le i \le N\}$ is the set of roots of a polynomial of 
degree $N\le n$
with  $C^{3n}$ coefficients.
By \cit!{7} there is an arrangement of these roots such 
that they become twice differentiable.
So claim \therosteritem{\nmb!{8}} follows.

The end of the proof is now similar to the end of the proof of
\therosteritem{D1}, where one uses claim \therosteritem{\nmb!{7}} (from the
proof of \therosteritem{D1}), claim \therosteritem{\nmb!{8}} instead of
\therosteritem{\nmb!{3}}, and the lemma above.
\qed\enddemo

\enddocument